\documentclass[reqno,11pt]{article}
\usepackage{mathrsfs}
\usepackage{amssymb}
\usepackage[usenames,dvipsnames]{xcolor}
\usepackage{amsthm}
\usepackage{amsmath}
\usepackage{amsfonts}
\usepackage{enumerate}
\usepackage{txfonts}
\usepackage{bm}
\usepackage{color}

\usepackage{graphicx}
\usepackage{tikz}
\usetikzlibrary{arrows,calc}
\tikzset{
	>=stealth',
	help lines/.style={dashed, thick},
	axis/.style={<->},
	important line/.style={thick},
	connection/.style={thick, dotted},}

\numberwithin{equation}{section}
\usepackage[unicode,breaklinks=true,colorlinks=true]{hyperref}
\newtheorem{theorem}{Theorem}[section]
\newtheorem{lemma}[theorem]{Lemma}
\newtheorem{corollary}[theorem]{Corollary}
\newtheorem{proposition}[theorem]{Proposition}

\theoremstyle{definition}
\newtheorem{definition}[theorem]{Definition}
\newtheorem{assumption}[theorem]{Assumption}
\newtheorem{example}[theorem]{Example}

\theoremstyle{remark}
\newtheorem{remark}[theorem]{Remark}

\def\Xint#1{\mathchoice
	{\XXint\displaystyle\textstyle{#1}}%
	{\XXint\textstyle\scriptstyle{#1}}%
	{\XXint\scriptstyle\scriptscriptstyle{#1}}%
	{\XXint\scriptscriptstyle\scriptscriptstyle{#1}}%
	\!\int}
\def\XXint#1#2#3{{\setbox0=\hbox{$#1{#2#3}{\int}$ }
		\vcenter{\hbox{$#2#3$ }}\kern-.6\wd0}}

\def\dashint{\Xint-}

\DeclareMathOperator*{\esssup}{ess\,sup}

\begin{document}
\title{On the existence of helical invariant solutions to steady Navier-Stokes equations}
\author{{Mikhail Korobkov\thanks{\noindent School of Mathematical Sciences,
Fudan University, Shanghai 200433, China; and Sobolev Institute of Mathematics, pr-t Ac. Koptyug, 4, Novosibirsk, 630090, Russia. E-mail: korob@math.nsc.ru.}}\and {Wenqi Lyu\thanks{Department of Mathematics,
	Southern University of Science and Technology, Shenzhen 518055, Guangdong Province, China. E-mail: lvwq@sustech.edu.cn.}}\and {Shangkun Weng\thanks{School of Mathematics and Statistics and Computational Science Hubei Key Laboratory, Wuhan University, Wuhan, Hubei Province, 430072, People's Republic of China. Email: skweng@whu.edu.cn.}}}
\date{\today}
\maketitle

\newcommand{\meas}{\mathop{\rm meas}}
\def\be{\begin{eqnarray}}
\def\ee{\end{eqnarray}}
\def\ba{\begin{aligned}}
\def\ea{\end{aligned}}
\def\bay{\begin{array}}
\def\eay{\end{array}}
\def\bca{\begin{cases}}
\def\eca{\end{cases}}
\def\p{\partial}
\newcommand{\dist}{\mathop{\rm dist}}
\renewcommand{\div}{\mathop{\rm div}}
\newcommand{\diam}{\mathop{\rm diam}}
\newcommand{\curl}{\mathop{\rm curl}}
\def\no{\nonumber}
\newcommand{\loc}{{\rm loc}}
\def\e{\varepsilon}
\def\ov{\boldsymbol{\omega}}
\def\de{\delta}
\def\De{\Delta}
\def\om{\omega}
\def\F{{\mathcal{F}}}
\def\we{\mathbf w}
\def\n{{\bf n}}
\def\ve{\mathbf w}
\def\eee{\mathbf e}
\def\fe{\mathbf f}
\def\ue{\mathbf u}
\def\Om{\Omega}
\def\f{\frac}
\def\Ha{{\mathfrak{H}}}
\def\N{{\mathbb{N}}}
\def\th{\theta}
\def\la{\lambda}
\def\lab{\label}
\def\w{{\mathbf{w}}}
\def\a{{\mathbf{a}}}
\def\h{{\mathbf{U}}}
\def\AA{{\mathbf{A}}}
\def\BB{{\mathbf{B}}}
\def\b{\bigg}
\def\var{\varphi}
\def\na{\nabla}
\def\Ti{{\mathscr{T}}}
\def\JI{{\mathscr{J}}}
\def\SI{{\mathscr{S}}}
\def\FI{{\mathscr{F}}}
\def\HI{{\mathscr{H}}}
\def\ka{\kappa}
\def\al{\alpha}
\def\La{\Lambda}
\def\ga{\gamma}
\def\const{{\rm const}}
\def\Ga{\Gamma}
\def\ti{\tilde}
\def\wti{\widetilde}
\def\wh{\widehat}
\def\ol{\overline}
\def\ul{\underline}
\def\Th{\Theta}
\def\R{\mathbb R}
\def\si{\sigma}
\def\wt{\widetilde}
\def\Si{\Sigma}
\def\oo{\infty}
\def\q{\quad}
\def\z{\zeta}
\def\co{:=}
\def\eqq{\eqqcolon}
\def\di{\displaystyle}
\def\bt{\begin{theorem}}
\def\et{\end{theorem}}
\def\bc{\begin{corollary}}
\def\ec{\end{corollary}}
\def\bl{\begin{lemma}}
\def\el{\end{lemma}}
\def\bp{\begin{proposition}}
\def\ep{\end{proposition}}
\def\br{\begin{remark}}
\def\er{\end{remark}}
\def\bd{\begin{definition}}
\def\ed{\end{definition}}
\def\bpf{\begin{proof}}
\def\epf{\end{proof}}
\def\bex{\begin{example}}
\def\eex{\end{example}}
\def\bq{\begin{question}}
\def\eq{\end{question}}
\def\bas{\begin{assumption}}
\def\eas{\end{assumption}}
\def\ber{\begin{exercise}}
\def\eer{\end{exercise}}
\def\mb{\mathbb}
\def\mbR{\mb{R}}
\def\mbZ{\mb{Z}}
\def\DD{\mb{D}}
\def\mc{\mathcal}
\def\mcS{\mc{S}}
\def\ms{\mathscr}
\def\lan{\langle}
\def\ran{\rangle}
\def\lb{\llbracket}
\def\rb{\rrbracket}
\begin{abstract}
In this paper, we investigate the nonhomogeneous boundary value problem for the steady Navier-Stokes equations in a helically symmetric spatial domain. When data is assumed to be helical invariant and satisfies the compatibility condition, we prove this problem has at least one helical invariant solution.
\end{abstract}



\section{Introduction}


Let $\Om=\Om_{0}\backslash\cup_{j=1}^{n}\bar{\Om}_{j}$ be a bounded multiply connected domain in $\mbR^3$ with $C^2$-smooth boundary $\p\Om= \bigcup_{j=0}^N \Ga_j$ consisting of $N+1$ disjoint components $\Ga_j=\p\Om_{j}, j=0,...,N$. Consider the nonhomogeneous boundary value problem for the steady Navier-Stokes equations
\be\lab{snsh}
\left\lbrace
\begin{split}
({\bf u}\cdot \na) {\bf u} + \na p =\Delta {\bf u} + {\bf f} \q & \text{in}\ \Om,\\
\text{div }{\bf u}=0 \q & \text{in}\ \Om,\\
{\bf u}={\bf a} \q & \text{on}\ \p\Om,
\end{split}
\right.
\ee
where ${\bf u}$ and $p$ are unknown velocity and pressure, ${\bf a}$ and {\bf f} are given boundary value and body force, we assume the viscous coefficient $\nu=1$ for simplicity. The boundary data should satisfy the compatibility condition
\be\lab{ch}
\int_{\p\Om} {\bf a}\cdot {\bf n} d S= \sum_{j=0}^{N} \int_{\Ga_j} {\bf a}\cdot {\bf n} d S= \sum_{j=0}^N \mc{F}_j=0,
\ee
where ${\bf n}$ is a unit outward normal vector to $\partial\Om$.

In his remarkable article $\cite{leray33}$, Leray proved the solvablity of $(\ref{snsh})$ when the flux of the velocity across each connected component $\Gamma_{j}$ of the boundary vanishes:
$$\int_{\Gamma_{j}} {\bf a\cdot n}\,dS=0, \quad j=0,...,N.$$

It remained an open problem as to whether the necessary condition $(\ref{ch})$ is sufficient for $(\ref{snsh})$ to be solvable. This is also called Leray's problem because it actually goes back to his paper $\cite{leray33}$.

Then, Leray's problem has been studied in many papers \cite{l1,l2,l3,l4,l5,l6,l7,l8,l9,l10,l11,l12,l13,l14}, only recently its solvability was proved in bounded 2D domains and for 3D axially symmetric case under the sole necessary condition $(\ref{ch})$ by M.Korobkov, K.Pileckas and R.Russo $\cite{kpr15annals}$. As far as we know, this problem for general 3D bounded domain remains open.

Given a positive number $\si$, we define the action of the helical group of transformations $G_{\sigma}$ on $\mbR^3$ by
\be\no
S_{\th,\si}(x)&=&\left(\bay{ll}
x_1 \cos\th + x_2 \sin \th\\
-x_1 \sin \th + x_2 \cos\th\\
\q x_3 +\f{\sigma}{2\pi} \th
\eay\right),\quad \th\in \mbR,
\ee
that is, a rotation around the $x_3$ axis with simultaneous translation along the $x_3$ axis. $G_{\si}$ is uniquely determined by $\si$, which we will call the {\it step}. We say that the smooth function $f(x)$ is helically symmetric or simply {\it helical}, if $f$ is invariant under the action of $G_{\si}$, i.e., $f(S_{\th,\si}(x))= f(x),\forall \theta\in\mbR$. Similarly, we say that the smooth vector field ${\bf u}(x)$ is helically symmetric, or simply helical, if it is covariant with respect to the action of $G_{\si}$, i.e., $M(\theta) {\bf u}(x)= {\bf u}(S_{\th,\si}(x))$ for all $\th \in \mbR$, where
\be\lab{p22}
M(\th)\co \left[\bay{lll}
\cos\th & \sin\th & 0\\
-\sin\th & \cos\th & 0\\
0 & 0 & 1
\eay\right].\ee

We learn from the definitions of helical function and vector field that they are $\sigma$-periodic in $x_{3}$ variable. A domain $\Omega\subset\mathbb{R}^{2}\times\mathbb{T}_{\sigma}$ is called a \textit{helical} domain if for each point $x\in\Omega$, $S_{\th,\si}(x)\in\Omega$ for any $\theta\in\mathbb{R}$. Here we denote by
$\mathbb{T}_{\sigma}=\mathbb{R}/\sigma\mathbb{Z}$ the corresponding 1-dimensional torus. In other words, the domain $\Omega$ is evolved from a two dimensional multiply connected domain $\mb{D}\subset \mbR^2$, hence
\be\no
\Om=\b\{\left(\begin{array}{ll}
	x_1\\
	x_2\\
	x_3
\end{array}\right): \text{$\exists (y_1,y_2)\in\mb{D}$ and $\theta\in\mbR$ such that } \left(\begin{array}{ll}
x_1\\
x_2\\
x_3
\end{array}\right)= \left(\begin{array}{ll}
y_1\cos\theta + y_2\sin\theta\\
-y_1\sin\theta + y_2 \cos\theta\\
\q\q \f{\si}{2\pi}\theta
\end{array}\right)\b\}.
\ee

There is an alternative definition of helical symmetry as follows. We rewrite a vector field ${\bf u}(x)= (u_1, u_2, u_3)(x_1,x_2,x_3)$ with respect to the moving orthonormal frame associated to standard cylindrical coordinates $(r,\theta,z)$,
\be\label{basis1}
{\bf e}_r= (\cos\theta, \sin\theta, 0),\q {\bf e}_{\th}= (-\sin\theta, \cos\theta, 0),\q {\bf e}_z= (0,0,1),
\ee
as ${\bf u}= u_r {\bf e}_r + u_{\th} {\bf e}_{\th} + u_z {\bf e}_z$, where $u_r, u_{\th}, u_z$ are functions of $(r,\theta, z)$. We introduce two new independent variables in place of $\theta$ and $z$:
\be\lab{p23}
\eta\co \f{\si}{2\pi}\theta +z,\quad \xi\co \f{\si}{2\pi}\theta-z.
\ee

A smooth function $p=p(r,\theta, z)$ is a helical function if and only if, when expressed in the $(r,\xi,\eta)$ variables, it is independent of $\xi$: $p=q(r,\f{\si}{2\pi}\theta+z)$, for some $q=q(r,\eta)$. Indeed, by definition $f$ is helical if and only if $f(S_{\rho,\sigma}(x))$ is actually independent of $\rho$, being equal to $f(x)$. Then by the Chain Rule $\f{d}{d\rho} f(S_{\rho,\sigma}(x))=0$ is equivalent to $\f{d}{d\xi} \ti{f}(r,\xi,\eta)=0$, where
\be\no
\ti{f}(r,\xi,\eta)= f\b(r\cos\b(\f{\pi}{\si}(\eta+\xi)\b), r\sin\b(\f{\pi}{\si}(\eta+\xi)\b), \f{\eta-\xi}{2}\b).
\ee

Similarly, a smooth vector field ${\bf u}$ is helical if and only if there exist $v_r, v_{\th}, v_z$ functions of $(r,\eta)$ such that $u_r(r,\theta,z)= v_r(r,\f{\si}{2\pi}\theta+z)$, $u_{\th}(r,\theta,z)= v_{\th}(r,\f{\si}{2\pi}\theta+z)$ and $u_z(r,\theta,z)= v_z(r,\f{\si}{2\pi}\theta+z)$. Since $u_r(r,\theta,z)= u_r (r,\theta+2\pi,z)$, then $v_r(r,\f{\si}{2\pi}\theta+z)= v_r(r,\f{\si}{2\pi}\theta+z + \si)$, that is to say, $v_r(r,\eta)$ is periodic in $\eta$ with period $\si$. Clearly, as $\sigma\to\infty$ a helical flow becomes 2D flow and when $\sigma=0$ it becomes axisymmetric flow (the coefficients $u_{r}$,$u_{\th}$,$u_{z}$ depend only on $r$ and $z$). Hence, a helical flow can be regarded as an interpolation between 2D and 3D axisymmetric flows.

It is well-known that the global regularity of unsteady Navier-Stokes equations is a longstanding open problem which is one of the seven Millennium Prize problems \cite{Feff}. But if the initial data is two-dimensional or axisymmetric without swirl ($u_{\theta}=0$), the global strong solution exists, see for example \cite{Tsai2018}. In general, the global well-posedness of Navier-Stokes equations with axisymmetric initial data is still open, see some partial results in \cite{Tsai2008,Sverak2009}. However, if the initial data is helical invariant, then the global regularity of the Leray-Hopf weak solution was proved by A.Mahalov, E.Titi and S.Leibovich \cite{mtl90arma}. In particular, a helical invariant function ${\bf u}$ satisfies the 2D Ladyzhenskaya inequality:
\begin{equation}
\|{\bf u}\|_{L^{4}(\Omega)}\leq C\sigma^{\frac{1}{4}}\|{\bf u}\|_{L^{2}(\Omega)}^{\frac{1}{2}}\|{\bf u}\|_{H^{1}(\Omega)}^{\frac{1}{2}},
\end{equation}
where $\Omega=\mathbb{R}^{2}\times\mathbb{T}_{\sigma}$. This coefficient $\sigma^{\frac{1}{4}}$ also reflects that the helical flow is an interpolation between 2D and 3D axisymmetric flows. There are more interesting iteratures concerning helical flows, see \cite{bfnn13siam,lmnnt14jdde,ettingertiti09siam,Sun2006,Liu2017,Jiu2017}.

This note is devoted to solve the Leray's problem in a helical domain with helical invariant data. More precisely, we prove the following theorem.

\bt\lab{main}{\sl
Assume that $\Omega\subset\mathbb{R}^{2}\times\mathbb{T}_{\sigma}$ is a helical domain with $C^2$ boundary. If ${\bf a}\in W^{3/2,2}(\partial\Omega)$ and ${\bf f}\in L^2(\Omega)$ are helical invariant functions, condition $(\ref{ch})$ is fulfilled, then problem $(\ref{snsh})$ admits at least one helical invariant weak solution.}
\et

\br\lab{reg}{\rm Under the hypothesis of Theorem~\ref{main}, every weak solution ${\bf u}$ is indeed a strong solution, i.e. ${\bf u}\in W^{2,2}(\Omega)$ (see, e.g., \cite{galdi11}). We need such regularity to apply Morse-Sard type theorem for regular level sets of Sobolev functions (see Theorem~\ref{MST1} and the commentary after the~formula~(\ref{boundary5})\,), so ${\bf a}\in W^{\frac{3}{2},2}(\partial\Omega)$ is necessary.}
\er

Our proof is motivated by M.Korobkov, K.Pileckas and R.Russo's approach which successfully solved the Leray's problem in 2D and 3D axisymmetric domain in \cite{kpr15annals}. It is plausible that we can extend this result to the helical case. The proof in \cite{kpr15annals} was based on Leray's contradiction arguments and
on the integration using the coarea formula along the level lines of the total head pressure.
Also they used  some consequence of Bernoulli's law for solutions~$\we$ to Euler equations with low regularity, namely, in \cite{kpr13,kpr15pisa}, the authors proved that for any compact connected set~$K\subset\Omega$
\begin{equation}\label{BI}
\biggl(\psi|_K\equiv\const\biggr)\Rightarrow \biggl(\Phi|_{K}\equiv\const\biggr),
\end{equation}
where $\psi$ is the corresponding stream function and $\Phi=p+\frac12|\we|^2$ is the total head pressure
(the identity for $\Phi$ is understood up to a negligible set of 1-dimensional measure zero).

The helical invariant functions form a closed subspace with respect to weak convergence, we can still use the contradiction arguments to derive a helical invariant weak solution $({\bf w},p)$ to the Euler equations
\be\lab{helicaleuler}
\left\lbrace
\begin{split}
	w_r\f{\p w_r}{\p r} + \left(\f{\si}{2\pi r}w_{\th}+ w_z\right) \f{\p w_r}{\p \eta} -\f{w_{\th}^2}{r}+ \f{\p p}{\p r} & =0,\\
	w_r\f{\p w_{\th}}{\p r}+ \left(\f{\si}{2\pi r}w_{\th}+ w_z\right) \f{\p w_{\th}}{\p \eta} +\f{w_r w_{\th}}{r} +\f{\si}{2\pi r}\f{\p p}{\p \eta} & =0,\\
	w_r\f{\p w_z}{\p r} + \left(\f{\si}{2\pi r}w_{\th}+ w_z\right) \f{\p w_z}{\p \eta} + \f{\p p}{\p \eta} & =0,\\
	\f{\p}{\p r}(r w_r) + \f{\p}{\p\eta}\left(\f{\si}{2\pi} w_{\th}+ r w_z\right) & =0,\\
      \we|_{\partial\Omega} & =0.
\end{split}
\right.
\ee
Unfortunately, although an equation
\be\lab{bernoulli}
\b(w_r\frac{\p}{\p r} + (\f{\si}{2\pi r}w_{\theta}+w_z)\frac{\p}{\p\eta}\b) \b(\f{1}{2}|{\bf w}|^2 + p\b)=0.
\ee
is still valid, the Bernoulli's law which appeared before is no longer true. To see this, let $w_{r}=0$, $w_{\th}=r$, $w_{z}=-\frac{\sigma}{2\pi}$, $p=\frac{r^{2}}{2}$. This is an explicit solution to equations $(\ref{helicaleuler})$, the stream function is constant since we always have $w_{r}=\frac{\sigma}{2\pi r}w_{\theta}+w_{z}=0$, but the total head pressure $\Phi=r^{2}+\frac{\sigma^{2}}{8\pi}$ is not a~constant. The similar problem appears in axially--symmetric case (see, e.g., \cite[page~746]{kpr18ma}\,), nevertheless, the Bernoulli identity~(\ref{BI})  will be satisfied here, due to zero boundary
conditions~$(\ref{helicaleuler})_{5}$ and since each straight line passing through a point of~$\Omega$ parallel to the axis of symmetry intersects~$\partial\Omega$. But this is not true anymore for helical domains, therefore, the Bernoulli identity~(\ref{BI})  fails in general
for the solutions to~(\ref{helicaleuler}) in helical case.

This is an obstacle in modifying the arguments in \cite{kpr15pisa,kpr15annals,kpr18ma}, since the Bernoulli's law is needed to separate the boundary components on which the total head pressure attains its supremum from the others, this step played a crucial role in solving Leray's problem in 3D axisymmetric domains.

On the other hand, if a function with a finite Dirichlet integral is helically symmetric, then its restriction to a two-dimensional hyperplane orthogonal to the axis of symmetry has a finite Dirichlet integral as well (see, e.g., the formula~(\ref{summ1})\,), i.e., the corresponding plane function has no singularities. This circumstance allows us to simplify the arguments and carry out the proof without using identities of the form~(\ref{BI}), which turned out to be so significant in the axisymmetric case (when the restriction of considered functions to the two-dimensional hyperplane have singularities near the~symmetry axis).

This paper is organized as follows. In section 2, we deal with the nonhomogeneous boundary values and present some properties of the Sobolev functions. In section 3, we first use Leray's contradiction arguments to derive a nontrivial solution to the helical Euler equations. Finally we construct a family of level sets of the total head pressure and deduce a contradiction via coarea formula.

\medskip

\noindent {\em Acknowledgment.} The authors are grateful to the referee for his/her careful review and valuable comments on the revision of the manuscript. The work of S.W. is partially supported by National Natural Science Foundation of China No.11701431, 11971307, 12071359.


\section{Preliminaries}


In this section, we first find a solenoidal and helical invariant extension of the nonhomogeneous boundary values. Then we review the Morse-Sard property of Sobolev functions.

By {\it a domain} we mean an open connected set. We use standard
notations for function spaces: $W^{k,q}(\Omega)$,
$W^{\alpha,q}(\partial\Omega)$, where $\alpha\in(0,1),
k\in{\mathbb N}_0, q\in[1,+\infty]$. In our notation we do not
distinguish function spaces for scalar and vector valued
functions; it is clear from the context whether we use scalar or
vector (or tensor) valued function spaces.

Denote by  $H(\Omega)$
the
closure of the set of all solenoidal smooth vector-functions having compact supports
in $\Omega$ with respect to the norm
$\|\we\|_{H(\Omega)}=\biggl(\int\limits_\Omega|\nabla\we|^2\biggr)^{\frac12}$.

\subsection{Extension of boundary values}

We need to use the following symmetry assumptions: \\

\begin{enumerate}
\item[$(H_{0})$] $\Omega\subset\mathbb{R}^3$ is a helical domain with $C^2$ boundary and $O_{x_3}$ is a symmetry axis of $\Omega$. \\

\item[$(H_{1})$] The assumptions $(H_{0})$ are fulfilled and both the boundary value ${\bf a}\in W^{3/2,2}(\partial\Omega)$ and ${\bf f}={\rm curl}\,{\bf b}\in L^2(\Omega)$ are helical invariant.
\end{enumerate}

\bl{\sl
If conditions $(H_{1})$ and $(\ref{ch})$ are fulfilled, then there exists a helical solenoidal extension ${\bf A}\in W^{2,2}(\Omega)$ of ${\bf a}$ with
\be\label{bound}
\|{\bf A}\|_{W^{2,2}(\Omega)}\leq c\|{\bf a}\|_{W^{3/2,2}(\partial\Omega)}.
\ee
}\el

\bpf
Let conditions $(H_{1})$ and $(\ref{ch})$ be fulfilled. Define by $P$ the hyperplane~$P=\{(y_1,y_2,0):y_1,y_2\in\R\}$. Then by construction
$$\DD=\Omega\cap P.$$
(cf. with the Introduction). Note, that by our assumptions, $\DD$ is a~bounded two dimensional domain with $C^2$-smooth boundary $\Si=\partial\DD=P\cap\partial\Omega$, \ $0\notin\overline{\DD}$, moreover,

{\sl every vector field $\tilde\AA:\DD\to\R^3$ can uniquely be extended to a helical invariant vector field ${\mathbf A}:\Omega\to\R^3$
such that ${\mathbf A}|_{\DD}=\tilde\AA$ \ on \ $\DD$ and } 
\be\label{bound}
\|{\bf A}\|_{W^{2,2}(\Omega)}\leq c\|{\tilde\AA}\|_{W^{2,2}(\DD)}.
\ee 
Indeed, if $\tilde\AA(r,\theta)=\tilde A_r(r,\theta)\eee_r+\tilde A_\theta(r,\theta)\eee_\theta+\tilde A_z(r,\theta)\eee_z$ with respect to the moving orthonormal frame associated to standard cylindrical coordinates $(r,\theta,z)$ \ (see~(\ref{basis1})\,), then $\AA$ can be defined as 
$\AA=A_r\eee_r+A_\theta\eee_\theta+A_z\eee_z$ with
\be\label{basis2}
A_r(r,\theta,z):=\tilde A_r\bigl(r,\theta+\frac{2\pi}{\sigma}z\bigr),\qquad
A_\theta(r,\theta,z):=\tilde A_\theta\bigl(r,\theta+\frac{2\pi}{\sigma}z\bigr),\qquad A_z(r,\theta,z):=\tilde A_z\bigl(r,\theta+\frac{2\pi}{\sigma}z\bigr)
\ee
(see the Introduction for the explanations). 
Using standard elementary formulas for divergence operator for plane vector fields in polar coordinate system $(r,\theta)$ on $\DD$ and for spatial vector fields in cylindrical coordinate system $(r,\theta,z)$ on $\Omega$, it is easy to check that for above helical extension procedure the following property holds:
 \be\label{basis3}
\biggl(\div\AA=0\mbox{\ \ in }\Omega\biggr)\Leftrightarrow
\biggl(\div\tilde\BB=0\mbox{\ \ in }\DD\biggr),
\ee
where we denote by $\tilde\BB$  the plane vector field $\tilde\BB=\tilde A_r\eee_r+\bigl(\tilde A_\theta+\frac{2\pi}{\sigma}\tilde A_z\bigr)\eee_\theta$ on $\DD$. 
Moreover, using classical elementary formulas for surface integrals, we have
 \be\label{basis4}\nonumber
\biggl(\int\limits_{\partial\Omega}{\bf A\cdot n}\,dS=0\biggr)\Leftrightarrow
\biggl(\int\limits_{\partial\DD}{\tilde\BB\cdot n}\,ds=0\biggr),
\ee
where $dS$ and $ds$ mean integration with respect to surface area and length respectively. 

Thus, the problem of a helical  solenoidal extension from the boundary of spatial domain~$\Omega$ can be reduced to the problem of solenoidal extensions for two dimensional vector fields in plane domains, which is well  known  classical fact (see, e.g., \cite{galdi11} Chapter III.3\,). This finishes  the proof of the~Lemma.
\epf

Then by the existence results of Stokes system \cite[Chapter 3]{galdi11}, we can find a unique solution ${\bf U}\in W^{2,2}(\Omega)$ to the Stokes problem such that ${\bf U-A}\in H(\Omega)\cap W^{2,2}(\Omega)$, and
\be\label{f}
\int_{\Omega}\nabla {\bf U}\cdot\nabla {\bm \varphi}dx=\int_{\Omega}{\bf f}\cdot{\bm \varphi}dx, \quad\forall{\bm \varphi}\in H(\Omega).
\ee

Moreover,
\be\label{u}
\|{\bf U}\|_{W^{2,2}(\Omega)}\leq c(\|{\bf a}\|_{W^{3/2,2}(\partial\Omega)}+\|{\bf f}\|_{L^2(\Omega)}).
\ee

Actually, uniqueness ensures the solution ${\bf U}$ is also a helical invariant function, since ${\bf A}$ is helical invariant. Indeed, for any $\xi_{0}\in[0,\sigma]$ denote ${\bf U}_{\xi_{0}}(r,\eta,\xi)={\bf U}(r,\eta,\xi-\xi_{0})$, then ${\bf U}_{\xi_{0}}-{\bf A}\in H(\Omega)\cap W^{2,2}(\Omega)$ since ${\bf A}$ is helical invariant. We learn from $(\ref{f})$
\be\label{f1}
\int_{\Omega}\nabla {\bf U}_{\xi_{0}}\cdot\nabla {\bm \varphi}dx=\int_{\Omega}\nabla {\bf U}\cdot\nabla {\bm \varphi}_{-\xi_{0}}dx=\int_{\Omega}{\bf f}\cdot{\bm \varphi}_{-\xi_{0}}dx=\int_{\Omega}{\bf f}_{\xi_{0}}\cdot{\bm \varphi}dx.
\ee
By our assumption, ${\bf f}$ is helical invariant, therefore,
${\bf f}_{\xi_{0}}={\bf f}$. We have ${\bf U}_{\xi_{0}}={\bf U}$, in another words, ${\bf U}$ is independent of $\xi$ and ${\bf U}$ is helical invariant.  Associated with the Stokes problem~\eqref{f}, there is a pressure term, which we denoted it to be $Q$.

Now we consider ${\bf w=u-U}, \tilde{p}=p-Q$, to simplify the notation, we still write $\tilde{p}$ as $p$, equation $(\ref{snsh})$ is equivalent to
\be
\begin{cases}
-\Delta {\bf w}+({\bf U}\cdot\nabla){\bf w}+({\bf w}\cdot\nabla){\bf w}+({\bf w}\cdot\nabla){\bf U}=-\nabla p-({\bf U}\cdot\nabla){\bf U}                   & \text{in} \ \Omega,\\
\text{div}\,{\bf w}=0                     & \text{in} \ \Omega,\\
{\bf w}=0                               & \text{on} \ \p\Omega.
\end{cases}
\ee

By a weak solution to problem $(\ref{snsh})$, we mean a function ${\bf u}\in W^{1,2}(\Omega)$ such that ${\bf w=u-U}\in H(\Omega)$ and for any ${\bm \varphi}\in H(\Omega)$,
\be\label{weak}
\begin{split}
\langle{\bf w},{\bm \varphi}\rangle_{H(\Omega)}= & -\int_{\Omega}({\bf U}\cdot\nabla){\bf U}\cdot{\bm \varphi}dx-\int_{\Omega}({\bf U}\cdot\nabla){\bf w}\cdot{\bm \varphi}dx \\
& -\int_{\Omega}({\bf w}\cdot\nabla){\bf w}\cdot{\bm \varphi}dx-\int_{\Omega}({\bf w}\cdot\nabla){\bf U}\cdot{\bm \varphi}dx.
\end{split}
\ee

By Riesz representation theorem, for any ${\bf w}\in H(\Omega)$ there exists a unique function $T{\bf w}\in H(\Omega)$ such that the right hand side of (\ref{weak}) is equivalent to $\langle T{\bf w},{\bm \varphi}\rangle_{H(\Omega)}$, for any ${\bm \varphi}\in H(\Omega)$. Moreover, $T$ is a compact operator by the well-known result in \cite{lady69}.

\begin{lemma}\label{th}
For any helical function ${\bf w}\in H(\Omega)$, $T{\bf w}\in H(\Omega)$ is also a helical function.
\end{lemma}
\bpf
The proof is the same as above. We omit the details here.
\epf

\subsection{Morse-Sard property of Sobolev functions}

The following lemma is concerned with the classical differentiablity properties of Sobolev functions.
\bl[\cite{dor}]
{\sl	If $\psi\in W^{2,1}(\mathbb{R}^{2})$, then $\psi$ is continuous and there exists a set $A_{\psi}$ such that $\mathfrak{H}^{1}(A_{\psi})=0$ and $\psi$ is differentiable (in the classical sense) at each $x\in\mathbb{R}\backslash A_{\psi}$. Furthermore, the classical derivative at these points $x$ coincides with $\nabla\psi(x)=\lim\limits_{r\to 0}\dashint_{B_{r}(x)}\nabla\psi(y)dy$, where $\lim\limits_{r\to 0}\dashint_{B_{r}(x)}|\nabla\psi(y)-\nabla\psi(x)|^{2}dy=0$.
}\el

Here and henceforth we denote by $\mathfrak{H}^1$ the
one-dimensional Hausdorff measure, i.e.,
$\mathfrak{H}^1(F)=\lim\limits_{t\to 0+}\mathfrak{H}^1_t(F)$,
where $\mathfrak{H}^1_t(F)=\inf\{\sum\limits_{i=1}^\infty {\rm
diam} F_i:\, {\rm diam} F_i\leq t, F\subset
\bigcup\limits_{i=1}^\infty F_i\}$.

The following Morse-Sard Theorem for Sobolev function has been proved by J.Bourgain, M.Korobkov and J.Kristensen \cite{bkk13,bkk15}.

\begin{theorem}\label{MST1}
	Let ${\DD}\subset\mathbb{R}^{2}$ be a bounded domain with Lipschitz boundary and $\psi\in W^{2,1}({\DD})$. Then:
	\begin{enumerate}
		\item[$(1)$] $\mathfrak{H}^{1}(\{\psi(x): x\in\bar{{\DD}}\backslash A_{\psi}\ \&\ \nabla\psi(x)=0\})=0$;
		\item[$(2)$] for every $\epsilon>0$ there exists $\delta>0$ such that for every set $U\subset\bar{{\DD}}$ with $\mathfrak{H}_{\infty}^{1}(U)<\delta$ the inequality $\mathfrak{H}^{1}(\psi(U))<\epsilon$ holds;
		\item[(3)] for $\mathfrak{H}^{1}$-almost all $y\in\psi(\bar{{\DD}})\subset\mathbb{R}$ the preimage $\psi^{-1}(y)$ is a finite disjoint family of $C^{1}$-curves $S_{j}$, $j=1,2,...,N(y)$. Each $S_{j}$ is either a cycle in ${\DD}$ (i.e., $S_{j}\subset{\DD}$ is homemorphic to the unit circle $\mathbb{S}^{1}$) or a simple arc with endpoints on $\partial{\DD}$ (in this case $S_{j}$ is transversal to $\partial{\DD}$).
	\end{enumerate}
\end{theorem}


\section{Proof of Theorem \ref{main}}

\subsection{Contradiction argument}


 Without loss of
generality, we may assume that ${\bf f}={\rm curl}\,{\bf b}\in
L^2(\Omega)$.\footnote{By the
Helmholtz--Weyl decomposition, $\fe$ can be represented as the sum
$\fe={\rm curl}\,{\bf b}+\nabla \varphi$ with ${\rm curl}\,{\bf
b}\in L^2(\Omega)$, and the
gradient part is included then into the pressure term (see, e.g.,
\cite{lady69}, \cite{galdi11}).} In particular,
\be\label{divv}
\div\fe=0\ee
in the sense of distributions.

It is obvious, (\ref{weak}) is equivalent to the operator equation ${\bf w}=T{\bf w}$ in $H(\Omega)$. By Leray-Schauder theorem, to prove the existence of weak solution is sufficient to show that the solutions of equation ${\bf w}^{(\lambda)}=\lambda T{\bf w}^{(\lambda)}$ are uniformly bounded with respect to $\lambda\in [0,1]$. In another word, the solution ${\bf w}\in H(\Omega)$ satisfies for any ${\bm \varphi}\in H(\Om)$,
\be\no
\begin{split}
\int_{\Omega}\nabla{\bf w}\cdot\nabla{\bm \varphi}dx =&-\lambda\int_{\Omega}({\bf U}\cdot\nabla){\bf U}\cdot{\bm \varphi}dx-\lambda\int_{\Omega}({\bf U}\cdot\nabla){\bf w}\cdot{\bm \varphi}dx \\
&-\lambda\int_{\Omega}({\bf w}\cdot\nabla){\bf w}\cdot{\bm \varphi}dx-\lambda\int_{\Omega}({\bf w}\cdot\nabla){\bf U}\cdot{\bm \varphi}dx
\end{split}
\ee
are uniformly bounded in $H(\Omega)$ with respect to $\lambda\in [0,1]$. \\

Assume this is false. Then there exist $\{\lambda_n\}_{n\in\mathbb{N}}\subset[0,1]$ and $\{{\bf \hat{w}}_n\}_{n\in\mathbb{N}}\subset H(\Omega)$ such that, for any ${\bm \varphi}\in H(\Omega)$,
\be\label{hat}
\begin{split}
\int_{\Omega}\nabla{\bf \hat{w}}_n\cdot\nabla{\bm \varphi}dx =& -\lambda_n\int_{\Omega}({\bf U}\cdot\nabla){\bf U}\cdot{\bm \varphi}dx-\lambda_n\int_{\Omega}({\bf U}\cdot\nabla){\bf \hat{w}}_n\cdot{\bm \varphi}dx \\
& -\lambda_n\int_{\Omega}({\bf \hat{w}}_n\cdot\nabla){\bf \hat{w}}_n\cdot{\bm \varphi}dx-\lambda_n\int_{\Omega}({\bf \hat{w}}_n\cdot\nabla){\bf U}\cdot{\bm \varphi}dx
\end{split}
\ee
and
\be\no
\lambda_n\to\lambda_{0}\in[0,1], J_n :=||{\bf \hat{w}}_n||_{H(\Omega)}\to\infty
\ee

After integration by part, $(\ref{hat})$ is equivalent to
\be\label{bypart}
\begin{split}
\int_{\Omega}\nabla{\bf \hat{w}}_n\cdot\nabla{\bm \varphi}dx=& \lambda_n\int_{\Omega}({\bf U}\cdot\nabla){\bm \varphi}\cdot{\bf U}dx+\lambda_n\int_{\Omega}({\bf U}\cdot\nabla){\bm \varphi}\cdot{\bf \hat{w}}_ndx \\
+& \lambda_n\int_{\Omega}({\bf \hat{w}}_n\cdot\nabla){\bm \varphi}\cdot{\bf \hat{w}}_ndx+\lambda_n\int_{\Omega}({\bf \hat{w}}_n\cdot\nabla){\bm \varphi}\cdot{\bf U}dx
\end{split}
\ee

Denote ${\bf w}_n={\bf \hat{w}}_n/J_n$. Since
\be\label{norm1}
\|{\bf w}_n\|_{H(\Omega)}=1,
\ee
 there exists a subsequence ${\bf w}_{n_l}$ converging weakly in $H(\Omega)$ to a function ${\bf w}\in H(\Omega)$. Then by compact embedding of Sobolev space, ${\bf w}_{n_l}$ converges strongly in $L^r(\Omega)$, for any $r\in [1,6)$. Taking ${\bm \varphi}=J_n^{-2}{\bf \hat{w}}_n$ in (\ref{bypart}), we get
\be\label{prelimit}
\int_{\Omega}|\nabla{\bf w}_n|^2 dx=\lambda_n\int_{\Omega}({\bf w}_n\cdot\nabla){\bf w}_n\cdot{\bf U}dx+J_n^{-1}\lambda_n\int_{\Omega}({\bf U}\cdot\nabla){\bf w}_n\cdot{\bf U}dx
\ee

Therefore, taking into account~(\ref{norm1}) and
passing to a limit as $n_l\to\infty$ in~(\ref{prelimit}), we obtain
\be\label{limit1}
1=\lambda_{0}\int_{\Omega}({\bf w}\cdot\nabla){\bf w}\cdot{\bf U}dx
\ee
in particular, this implies $\lambda_0\in (0,1]$.

Let us return to $(\ref{bypart})$, for any ${\bm \zeta}\in W_{0}^{1,2}(\Omega)$, consider the linear functional
\be\label{rn1}
\begin{split}
R_n({\bm \zeta})=\int_{\Omega}\nabla{\bf \hat{w}}_n\cdot\nabla{\bm \zeta}dx & -\lambda_n\int_{\Omega}({\bf U}\cdot\nabla){\bm \zeta}\cdot{\bf U}dx-\lambda_n\int_{\Omega}({\bf U}\cdot\nabla){\bm \zeta}\cdot{\bf \hat{w}}_ndx \\
& -\lambda_n\int_{\Omega}({\bf \hat{w}}_n\cdot\nabla){\bm \zeta}\cdot{\bf \hat{w}}_ndx-\lambda_n\int_{\Omega}({\bf \hat{w}}_n\cdot\nabla){\bm \zeta}\cdot{\bf U}dx
\end{split}
\ee

It is obvious that $R_n({\bm \zeta})$ has the following estimate
\be\label{functional}
|R_n({\bm \zeta})|\leq C(\|{\bf \hat{w}}_n\|_{H(\Omega)}+\|{\bf \hat{w}}_n\|_{H(\Omega)}^2+\|{\bf U}\|_{W^{1,2}(\Omega))^2}\|{\bm \zeta}\|_{H(\Omega)}
\ee

Again $(\ref{bypart})$ implies
\be\no
R_n({\bm \varphi})=0,\quad\forall{\bm \varphi}\in H(\Omega).
\ee

Therefore, there exists a unique function $\hat{p}_n\in \hat{L}^2(\Omega)=\{q\in L^2(\Omega):\int_{\Omega}q(x)dx=0\}$ such that
\be\label{rn2}
R_n({\bm \varphi})=\int_{\Omega}\hat{p}_n\text{div}\,{\bm \varphi}dx,\quad\forall{\bm \varphi}\in H(\Omega).
\ee
actually $\hat{p}_{n}$ is a helical function, this can be established by uniqueness since the functions ${\bf \hat{w}}_n$ and ${\bf U}$ on the RHS of $(\ref{rn1})$ are helical invariant. The proof of this part is similar to $(\ref{f1})$. Then we can learn from $(\ref{functional})$,
\be\label{pnorm}
\begin{split}
\|\hat{p}_n\|_{L^2(\Omega)} & \leq C\|R_n\|_{H^{-1}(\Omega)} \\
& \leq C(\|{\bf \hat{w}}_n\|_{H(\Omega)}+\|{\bf \hat{w}}_n\|_{H(\Omega)}^2+\|{\bf a}\|_{W^{2,3/2}(\p\Omega)}^2+\|{\bf f}\|_{L^2(\Omega)}^2).
\end{split}
\ee
where we have used $(\ref{u})$.

Combine $(\ref{rn1})-(\ref{rn2})$, the pair $({\bf \hat{w}}_n,\hat{p}_n)$ satisfies for any ${\bm \zeta}\in W_{0}^{1,2}(\Omega)$
\be\no
\int_{\Omega}\nabla{\bf \hat{w}}_n\cdot\nabla{\bm \zeta}dx-\lambda_n\int_{\Omega}({\bf U}\cdot\nabla){\bm \zeta}\cdot{\bf U}dx-\lambda_n\int_{\Omega}({\bf U}\cdot\nabla){\bm \zeta}\cdot{\bf \hat{w}}_ndx \\
\label{hatp}-\lambda_n\int_{\Omega}({\bf \hat{w}}_n\cdot\nabla){\bm \zeta}\cdot{\bf \hat{w}}_ndx-\lambda_n\int_{\Omega}({\bf \hat{w}}_n\cdot\nabla){\bm \zeta}\cdot{\bf U}dx=\int_{\Omega}\hat{p}_n\text{div}\,{\bm \zeta}dx
\ee

Let ${\bf \hat{u}}_n={\bf \hat{w}}_n+{\bf U}$, then $(\ref{hatp})$ is equivalent to for any ${\bm \zeta}\in W_{0}^{1,2}(\Omega)$
\be\no
\int_{\Omega}\nabla{\bf \hat{u}}_n\cdot\nabla{\bm \zeta}dx-\int_{\Omega}\hat{p}_n\text{div}\,{\bm{\zeta}}dx=-\lambda_n\int_{\Omega}({\bf \hat{u}}_n\cdot\nabla){\bf \hat{u}}_n\cdot{\bm \zeta}dx+\int_{\Omega}{\bf f}\cdot{\bm \zeta}dx
\ee
where we have used $(\ref{f})$.

In another word, the pair $({\bf \hat{u}}_n,\hat{p}_n)$ satisfies the following Navier-Stokes equations
\be\label{hatns}
\begin{cases}
	-\Delta {\bf \hat{u}}_n+\nabla \hat{p}_n=-\lambda_n({\bf \hat{u}}_n\cdot\nabla){\bf \hat{u}}_n+{\bf f}                   & \text{in} \ \Omega,\\
	\text{div}\,{\bf \hat{u}}_n=0              & \text{in} \ \Omega,\\
	{\bf \hat{u}}_n={\bf a}                        & \text{on} \ \p\Omega.
\end{cases}
\ee

By the well-known $L^q$-estimates of Stokes systems (see, e.g.,
\cite[P.\,279,\,Theorem\,IV.4.1]{galdi11}\,), equations $(\ref{hatns})$ imply
\be\label{stokes}
\begin{split}
\|{\bf \hat{u}}_n\|_{W^{2,3/2}(\Omega)}+\|\hat{p}_n\|_{W^{1,3/2}(\Omega)} & \leq C(\|({\bf \hat{u}}_n\cdot\nabla){\bf \hat{u}}_n\|_{L^{3/2}(\Omega)}+\|{\bf f}\|_{L^{2}(\Omega)}+\|{\bf a}\|_{W^{3/2,2}(\p\Omega)}) \\
& \leq C(\|{\bf \hat{u}}_n\|_{W^{1,2}(\Omega)}^2+\|{\bf f}\|_{L^{2}(\Omega)}+\|{\bf a}\|_{W^{3/2,2}(\p\Omega)}) \\
& \leq C(\|{\bf \hat{w}}_{n}\|_{H(\Omega)}^2+\|{\bf a}\|_{W^{3/2,2}(\p\Omega)}^2+\|{\bf f}\|_{L^{2}(\Omega)}+\|{\bf a}\|_{W^{3/2,2}(\p\Omega)})
\end{split}
\ee
note that we have used $(\ref{pnorm})$ in the last step.



Let ${\bf u}_n=J_n^{-1}{\bf \hat{u}}_n$ and $p_n=\lambda_{n}^{-1}J_n^{-2}\hat{p}_n$. It is obvious that $\|{\bf u}_n\|_{W^{1,2}(\Omega)}$ and $\|p_n\|_{W^{1,3/2}(\Omega)}$ are uniformly bounded, and
\be\label{preeuler}
\begin{cases}
	-\nu_n\Delta {\bf u}_n+({\bf u}_n\cdot\nabla){\bf u}_n+\nabla p_n={\bf f}_n                   & \text{in} \ \Omega,\\
	\text{div}\,{\bf u}_n=0              & \text{in} \ \Omega,\\
	{\bf u}_n={\bf a}_n                        & \text{on} \ \p\Omega.
\end{cases}
\ee
where $\nu_n=\lambda_{n}^{-1}J_n^{-1}$, ${\bf f}_n=\lambda_{n}^{-1}J_n^{-2}{\bf f}$, and ${\bf a}_n=J_n^{-1}{\bf a}$. By the observation above, we can extract weakly convergent subsequences ${\bf u}_n\rightharpoonup {\bf w}$ in $W^{1,2}(\Omega)$ and $p_n\rightharpoonup p$ in $W_{loc}^{1,3/2}(\Omega)\cap L^2(\Omega)$. The pair $({\bf w}, p)$ satisfies the following Euler equations
\be\label{euler}
\begin{cases}
	({\bf w}\cdot\nabla){\bf w}+\nabla p=0    & \text{in} \ \Omega,\\
	\text{div}\,{\bf w}=0              & \text{in} \ \Omega,\\
	{\bf w}=0                        & \text{on} \ \p\Omega.
\end{cases}
\ee

We summarize the above results as follows.

\bl \label{NSkk2}{\sl
Assume that $\Om\subset\mathbb{R}^{2}\times\mathbb{T}_{\sigma}$ is a helical domain with $C^2$ boundary $\p\Om$, ${\bf f}={\rm curl}\,{\bf b}$, ${\bf b}\in W^{1,2}(\Om)$ and ${\bf a}\in W^{3/2,2}(\p\Om)$ are both helical functions and $(\ref{ch})$ is valid. If the assertion of Theorem $\ref{main}$ is false, then there exist ${\bf w},p$ with the following properties:
\begin{enumerate}
\item[${\bf (E\textendash H)}$] The helical functions ${\bf w}\in W_{0}^{1,2}(\Om)$, $p\in W^{1,3/2}(\Om)$ satisfy the Euler equations $(\ref{euler})$, moreover equation $(\ref{limit1})$ holds.

\item[${\bf(E\textendash NS\textendash H)}$] Conditions ${\bf (E\textendash H)}$ are satisfied, and there exist sequences of helical functions ${\bf u}_n\in W^{1,2}(\Om)$, $p_n\in W^{1,3/2}(\Om)$ and numbers $\nu_n\to 0^{+}$, $\lambda_n\to \lambda_0 >0$ such that the norms $\|{\bf u}_n\|_{W^{1,2}(\Om)}$, $\|p_n\|_{W^{1,3/2}(\Om)}$ are uniformly bounded, the pair $({\bf u}_n, p_n)$ satisfies $(\ref{preeuler})$, and

$\|\nabla {\bf u}_n\|_{L^2(\Om)}\to 1$, \ \ ${\bf u}_n\rightharpoonup {\bf w}$ in\ $W^{1,2}(\Om)$, \ \  $p_n\rightharpoonup p$ in $W^{1,3/2}(\Omega)$.
\end{enumerate}
Moreover,
\be\label{reg-ass1}
{\bf u}_{n} \in W^{2,2}(\Om)\qquad\mbox{ and }\quad p_n\in W_{loc}^{2,1}(\Om).
\ee}
\el
Note, that the last inclusion follows from~(\ref{divv}) and from  the fact that locally the laplacian $\Delta p$ belongs to the Hardy space (see, e.g.,~\cite[Section~5.2]{kpr13} for more detailed description).


\subsection{Euler equations}


In this subsection we discuss the properties satisfied by the helical invariant solution to Euler equation $(\ref{euler})$.

Define by $P$ the hyperplane~$P=\{(x_1,x_2,0):x_1,x_2\in\R\}$. Then by construction
$$\DD=\Omega\cap P.$$

It follows from the helical symmetry of the vector field $\we$ that $\|\we\|_{L^q(\Omega)}=\sigma\|\we(\cdot,0)\|_{L^q(\DD)}$, hence we have
\begin{equation}
\label{summ1}\nabla\we\in L^2(\DD), \qquad\we\in L^q(\DD) \quad \forall q<\infty,
\end{equation}
consequently, from Euler system we have
\begin{equation}
\label{summ2}\nabla p\in L^q(\DD)  \quad \forall q<2.
\end{equation}


In this section and below for any set $S\subset\mathbb{R}^2$ we define $\widetilde{S}\subset\mathbb{R}^2\times\mathbb{T}_{\sigma}$ to be the three dimensional set which is evolved from $S$, i.e
\be\label{evvv}
\widetilde{S}=  \left\{\left(\begin{array}{ll}
y_1\cos\theta + y_2\sin\theta\\
-y_1\sin\theta + y_2 \cos\theta\\
\q\q \f{\si}{2\pi}\theta
\end{array}\right): (y_1,y_2)\in S,\ \theta\in\R\right\}.
\ee

Denote $\Si_j:=P\cap \Gamma_j$. Clearly, $\widetilde{\Sigma}_{j}=\Gamma_{j}$, $j=0,\cdots,N$, and
\be\label{bound1}
\partial \DD=\bigcup\limits_{j=0}^N \Si_j.
\ee

The next statement  was proved in \cite[Lemma 4]{KaPi1} and in
\cite[Theorem 2.2]{aiumj84}.

\bl\label{pconst}{\sl
If ${\bf (E\textendash H)}$ are satisfied, then

$$\exists\check{p}_{j}\in\mathbb{R}: p(x)\equiv \check{p}_{j} \ \ \ \mbox{for\ $\Ha^2$-almost\ all \ }x\in\Gamma_{j},\ j=0,...,N.$$

In particular, by helical symmetry,
$$p(x)\equiv \check{p}_j \ \ \ \mbox{for\ $\Ha^1$-almost\ all \ }x\in{\Sigma}_j,\ j=0,...,N.$$}
\el

By simple calculation from $(\ref{limit1})$, $(\ref{euler})_1$ and Lemma $\ref{pconst}$, it follows that

\bc{\sl
If conditions ${\bf(E\textendash NS\textendash H)}$ are satisfied, then
\be\label{plinearc}
-\frac{1}{\lambda_0}=\sum_{j=0}^{N}\check{p}_j\int_{\Gamma_j}{\bf a\cdot n}dS=\sum_{j=0}^{N}\check{p}_j\mc{F}_{j}.
\ee
}\ec

Set $\Phi_n=p_n+\frac{1}{2}|{\bf u}_n|^2$ and $\Phi=p+\frac{1}{2}|{\bf w}|^2$. By the properties of Sobolev functions best representatives for ${\bf w}$, $\Phi$ (see \cite{ep15}), we get the following.

\bl\lab{sobolev}
{\sl If conditions ${\bf (E\textendash H)}$ hold, then there exists a set $A_{{\bf w}}\subset\DD$ such that
\begin{enumerate}
	\item[(1)]$\Ha^1(A_{{\bf w}})=0$;
	\item[(2)]for all $x\in \DD\backslash A_{\bf w}$,
	\be
	\lim\limits_{\rho\to 0}\ \dashint_{B_{\rho}(x)}|{\bf w}(y)-{\bf w}(x)|^2dy=\lim\limits_{\rho\to 0}\ \dashint_{B_{\rho}(x)}|\Phi(y)-\Phi(x)|^2dy=0;
	\ee
	\item[(3)] for every $\epsilon>0$, there exists an~open set set $U\subset \R^2$ with $\Ha_{\infty}^{1}(U)<\epsilon$, $A_{\bf w}\subset U$ and such that the functions ${\bf w}$, $\Phi$ are continuous on $\bar{\DD}\backslash U$.
\end{enumerate}
}\el


\subsection{Obtaining a contradiction} We consider two possible cases.


\begin{enumerate}

\item[$(a)$] The maximum of $\Phi$ is attained on the boundary $\p\Om$:
\be\label{a}
\max_{j=0,...,N}\check{p}_{j}=\esssup\limits_{x\in\Omega}\Phi(x).
\ee

\item[$(b)$] The maximum of $\Phi$ is not attained on the boundary $\p\Om$\footnote{$\esssup\Phi=\oo$ is not excluded.}:
\be\label{b}
\max_{j=0,...,N}\check{p}_{j}<\esssup\limits_{x\in\Omega}\Phi(x).
\ee

\end{enumerate}

\subsubsection{If $(a)$ happens.}
Add to the pressure a constant such that $\esssup\limits_{x\in\Omega}\Phi(x)=0$.

Without loss of generality, we can renumerate the boundary components
in such a way that
\begin{equation}\label{phi=0}
\check{p}_0=\check{p}_1=\dots=\check{p}_M=0,\qquad 0\le M<N,
\end{equation}
\begin{equation}\label{phi-0}
\check{p}_j<0,\quad j=M+1,\dots,N.
\end{equation}
Denote $p_*=\frac13\max\limits_{j=M+1,...,N}\check{p}_j$. Then by construction
\begin{equation}\label{phi-00}
0>p_*>\check{p}_j,\qquad \forall j=M+1,\dots,N.
\end{equation}

For sufficiently small parameter $h>0$ and $j\in\{0,\dots,N\}$ denote
$\Si_{jh}=\{x\in\DD:\dist(x,\Si_j)=h\}$, \ $\DD_{jh}=\{x\in\DD:\dist(x,\Si_j)<h\}$, and
\begin{equation}\label{boundary00}
\Si=\Si_0\cup\dots\cup\Si_M,\qquad \Si_h=\bigcup\limits_{j=0}^M\Si_{jh}=\bigl\{x\in\DD:\dist(x,\Si)=h\bigr\},\qquad\DD_h=\bigcup\limits_{j=0}^M\DD_{jh}=\bigl\{x\in\DD:\dist(x,\Si)<h\bigr\},
\end{equation}
$$\Si^-=\Si_{M+1}\cup\dots\cup\Si_N,\qquad\Si^-_h=\bigcup\limits_{j=M+1}^N\Si_{jh}=\bigl\{x\in\DD:\dist(x,\Si^-)=h\bigr\},$$
$$\DD^-_h=\bigcup\limits_{j=M+1}^N\DD_{jh}=\bigl\{x\in\DD:\dist(x,\Si^-)<h\bigr\}.$$
In particular, we have
\begin{equation}\label{hh00}
\bigl\{x\in\DD:\dist(x,\partial\DD)=h\bigr\}=\Si_h\cup\Si^-_h.
\end{equation}
Since the distance function $\dist(x,\partial\DD)$ is
$C^1$--regular and the norm of its gradient is equal to one in the
neighborhood of $\partial\Omega$, there is a constant $\delta_0>0$
such that for every positive $h\le\delta_0$ the set $\Si_{jh}$ is a~$C^1$-smooth  curve homeomorphic to the~circle, below we will call such curves as {\it cycles}. Respectively,
${\Si}_h$ is a disjoint union of cycles.

By direct calculations, (\ref{euler}) implies
\begin{equation}\label{grthpax}\nabla\Phi=\ve\times\ov\quad\mbox{ in }\Omega,
\end{equation}
where $\ov=\curl\ve$, i.e.,
$$\ov=(\omega_r,\omega_\theta,\omega_z)=\bigl(-\frac{\partial v_\theta}{\partial z},\ \frac{\partial v_r}{\partial z}-\frac{\partial v_z}{\partial r},\ \frac{v_\theta}r+\frac{\partial v_\theta}{\partial r}\bigr).$$
Set $\omega(x)=|\ov(x)|$.

From conditions ${\bf(E\textendash NS\textendash H)}$ and from helical symmetry by construction we obtain
\begin{equation}\label{wwcc1}{\bf u}_n\rightharpoonup {\bf w}\quad\mbox{ in }W^{1,2}(\DD),
\end{equation}
\begin{equation}\label{wwcc2}p_n\rightharpoonup p\quad\mbox{ in }W^{1,3/2}(\DD).
\end{equation}

Then from Theorem 3.2 in \cite{aiumj84}	 (see also \cite[Lemma 3.3]{kpr13}) we obtain
\bl
\label{regc-ax} {\sl There exists
a~subsequence~$\Phi_{k_l}$ such that  $\Phi_{k_l}|_{\Si_{jh}}$ converges to
$\Phi|_{\Si_{jh}}$ uniformly for almost all $h\in(0,\delta_0)$ and for all $j=\{0,1,\dots,N\}$.}
\el

Below we assume (without loss of generality) that the subsequence $\Phi_{k_l}$ coincides with the whole sequence $\Phi_{k}$.

The value~$h\in(0,\delta_0)$ will be called
{\it regular}, if it satisfies the assertion of Lemma~\ref{regc-ax}, i.e., if
\begin{equation}\label{ucc1}\Phi_{k}|_{\Si_{jh}}\rightrightarrows \Phi|_{\Si_{jh}}\qquad\forall j=\{0,1,\dots,N\}.
\end{equation}

\bl
\label{reg-2} {\sl There exists a measurable set $\HI\subset(0,\delta_0)$
such that
\begin{enumerate}[(i)]
\item \ each value $h\in\HI$ is regular;
\item \ density of $\HI$ at zero equals~1:
\begin{equation}\label{md1} \frac{\meas\bigl(\HI\cap[0,h]\,\bigr)}{h}\to 1\qquad\mbox{ as \ }h\to0+;
\end{equation}
\item \qquad $\varlimsup\limits_{\HI\ni h\to0+}\sup\limits_{x\in \Si_{jh}}\bigl|\Phi(x)-\check{p}_j\bigr|=0\qquad\forall j=\{0,1,\dots,N\}.$
\end{enumerate}
}
\el

\bpf Fix $j\in\{0,1,\dots,N\}.$ \ Since $\we|_{\partial\DD}\equiv 0$ and $\nabla\we\in L^2(\DD)$, by Hardy inequality
we have
$$\int\limits_{\DD_{jh}}|\we|^2=o(h^2).$$
Then by H\"{o}lder inequality,
$$\int\limits_{\DD_{jh}}|\we|\cdot|\nabla\we|=o(h).$$
From Euler equations we have
$$\int\limits_{\DD_{jh}}|\nabla\Phi|\le 2\int\limits_{\DD_h}|\we|\cdot|\nabla\we|=o(h).$$
Using Fubini theorems, we can rewrite the last estimate as
$$\int\limits_{\DD_{jh}}|\nabla\Phi|=\int\limits_0^h\biggl(\int\limits_{\Si_{jt}}|\nabla\Phi|\,ds\biggr)dt=o(h).$$
The last estimate implies easily that there exists a measurable set $\HI\subset(0,\delta_0)$ such that
the density of $\HI$ at zero equals~1 (e.g., (\ref{md1}) holds\,), the restriction $\Phi|_{\Si_{jh}}$ is an absolute continuous function of one variable for all~$h\in\HI$,  and
 \begin{equation}\label{md-a1}\varlimsup\limits_{\HI\ni h\to0+}\int\limits_{\Si_{jh}}|\nabla\Phi|\,ds=0.
\end{equation}
Denote by $\overline{\Phi}_h$ the mean value of $\Phi$ over the curve $\Si_{jh}$.
Then from~(\ref{md-a1}) we obtain immediately that
 \begin{equation}\label{md-a2}\varlimsup\limits_{\HI\ni h\to0+}\max\limits_{x\in\Si_{jh}}|\Phi(x)-\overline{\Phi}_h|=0.
\end{equation}

Since by our assumptions $\partial\DD$ is $C^2$-smooth, the curves $\Sigma_{ih}$ are uniformly
$C^1$-smooth for all sufficiently small $h<\delta$. Then by  well-known classical Sobolev theorems,
the continuous trace operator
$$T_{jh}:W^{1,3/2}(\DD)\ni f\mapsto f|_{\Sigma_{jh}}\in L^2(\Sigma_{jh})$$
is well defined. Moreover, for any $f\in W^{1,3/2}(\DD)$  the real-valued function
 \begin{equation}\label{md-a3}[0,\delta]\ni h\mapsto \int\limits_{\Sigma_{jh}}|f|\,ds
\end{equation}
is continuous.

Recall, that by Lemma~\ref{pconst}  the trace identity $\Phi|_{\Si_j}\equiv\check{p}_j$ holds. Since pressure is defined up to an~additive constant, we can assume, without loss of generality, that $\check{p}_j=0$, i.e.,
 \begin{equation}\label{md-a4}\Phi|_{\Si_j}=0
\end{equation} in the sense of trace of Sobolev spaces. Hence, from (\ref{md-a3})--(\ref{md-a4}) we obtain
\begin{equation}\label{md-a44}\int\limits_{\Sigma_{jh}}|\Phi|\,ds\to0\qquad\mbox{ as }\ h\to0+.
\end{equation}
This implies $\overline{\Phi}_h\to0$ as $h\to0+$. The last convergence together with the formula~(\ref{md-a1})
give us the desired asymptotic
\begin{equation}\label{md-a7}\varlimsup\limits_{\HI\ni h\to0+}\max\limits_{x\in\Si_{jh}}|\Phi(x)|=0.
\end{equation}
The Lemma is proved.\epf

By Lemmas~\ref{regc-ax}--\ref{reg-2}, decreasing $\HI$ if necessary and taking sufficiently large~$k$, we can assume, without loss of generality, that

\bc\label{cest1}{\sl For all $h\in\HI$ we have
\be\label{fest1}\Phi_k(x)<p_*\qquad\forall x\in \Si^-_{h}.
\ee}\ec

Also, by Lemma~\ref{reg-2} we obtain
\bc\label{cest2}{\sl For all $\e>0$ there exists $\delta_\e>0$ and $k_\e\in\N$ such that $\forall h\in\HI\cap(0,\delta_\e)$ and $\forall k\ge k_\e$ we have
\be\label{fest2}\Phi_k(x)>-\e\qquad\forall x\in \Si_{h}.
\ee}\ec

\

Denote $t_*=-p_*$. Let $0<t_0<t_*$. The next geometrical object
plays an~important role in the estimates below: for
$t\in(t_0,t_*)$ and for all sufficiently large~$k$ we define the level set $S_k(t,t_0)\subset
\{x\in\DD:\Phi_k(x)=-t\}$ separating boundary components  $\Si$ from $\Si^-$ as follows. Namely, take
 $\epsilon=\frac12t_0$ and the corresponding parameters $\delta_\epsilon>0$, $k_\circ=k_\epsilon\in\N$  such that
(\ref{fest1})--(\ref{fest2}) holds for all $k\ge  k_\circ$ and  $h\in\HI\cap(0,\delta_\epsilon)$. Fix a number
${h_\circ}\in\HI\cap(0,\delta_\epsilon)$.
In particular, we have
\begin{equation}\label{boundary2}
\forall t\in(t_0,t_*)\ \ \forall k\ge k_\circ\quad
\biggl(\Phi_k|_{\Si_{h_\circ}}> -t,\quad\Phi_k|_{\Si^-_{h_\circ}}< -t\biggr).
\end{equation}

For $k\ge  k_\circ$, $j=0,1,\dots,M$, and $t\in(t_0,t_*)$
\begin{equation}\label{deff1}
\begin{array}{lcr}
\mbox{\it denote by
$W^j_k(t_0;t)$ the connected component of the open set $\bigl\{x\in
\DD:\dist (x,\partial\DD)> h_\circ\  \& \ \Phi_k(x)>-t\bigr\}$}
\\
\mbox{\it \qquad\qquad \qquad\qquad such that
$\partial W^j_k(t_0;t)\supset \Si_{j{h_\circ}}$ \ \ \  (see Fig.1), }
\end{array}
\end{equation}
and put
$$W_k(t_0;t)=\bigcup\limits_{j=0}^M W^j_k(t_0;t),\qquad
S_k(t_0;t)=\bigl(\partial W_k(t_0;t)\bigr)\setminus \Si_{h_\circ}.$$
By construction, (see Fig.1),
\begin{equation}\label{boundary6}\partial W_k(t_0;t)=\Si_{h_\circ}\cup S_k(t_0;t)
\end{equation}  and
\begin{equation}\label{boundary3}
\partial W_k(t_0;t)\subset\bigl\{x\in\DD:\dist(x,\partial\DD)> {h_\circ}\  \& \ \Phi_k(x)=-t\bigr\}\cup
\bigl\{x\in\DD:\dist(x,\partial\DD)={h_\circ}\  \& \ \Phi_k(x)\ge-t\bigr\}.
\end{equation}
\begin{center}
\includegraphics[height=6cm, width=7.5cm]{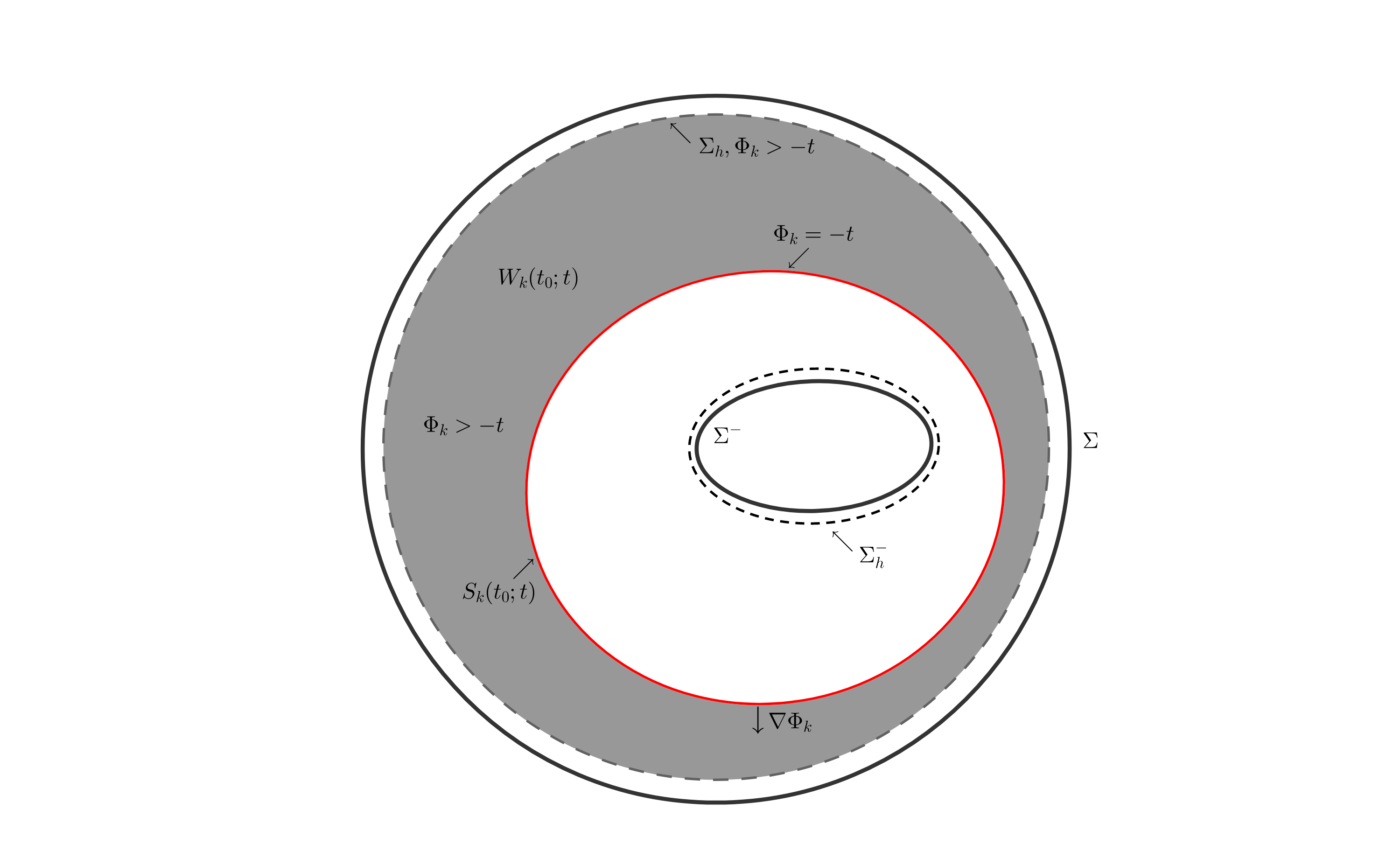}\\
{\small Figure 1: $M=0, N=1$.}
\end{center}

Therefore, by definition of $S_k(t_0;t)$ and in virtue of the identity~(\ref{hh00}),
\begin{equation}\label{boundary4}
S_k(t_0;t)\subset\bigl\{x\in\DD:\dist(x,\partial\DD)> {h_\circ}\  \& \ \Phi_k(x)=-t\bigr\}\cup \bigl\{x\in\Si^-_{h_\circ}:\Phi_k(x)\ge-t\bigr\}.
\end{equation}
But the last set $\{x\in\Si^-_{h_\circ}:\Phi_k(x)\ge -t\}$ is empty because of~(\ref{boundary2}). Therefore,
\begin{equation}\label{boundary5}
S_k(t_0;t)\subset\bigl\{x\in\DD:\dist(x,\partial\DD)>{h_\circ}\  \& \ \Phi_k(x)=-t\bigr\}.
\end{equation}

Since by ${\bf(E\textendash NS\textendash H)}$ each $\Phi_k$ belongs to   $W^{2,1}_{\loc}(\DD)$,
by the Morse-Sard theorem for Sobolev functions (see assertion~(iii) of Theorem~\ref{MST1}\,) we have that for almost all $t\in(t_0,t_*)$
the level set $S_k(t_0;t)$ consists of finitely many
$C^1$-cycles and $\Phi_k$ is differentiable (in classical sense)
at every point~$x\in S_k(t_0;t)$ with $\nabla\Phi_k(x)\ne0$.
The values $t\in(t_0,t_*)$ having the above property will be called
$k$-{\it regular}. (Note that $W_k(t_0;t)$ and
$S_k(t_0;t)$ are well defined for all $t\in(t_0,t_*)$ and $k\ge
 k_\circ= k_\circ(t_0)$.)

Recall that for a set $A\subset \DD$ we denote by $\widetilde A$
the three dimensional set in $\mathbb{T}_{\sigma}$ which is evolved from $A$ (see (\ref{evvv})\,).
By construction, for every regular value~$t\in(t_0,t_*)$ the set
$\widetilde S_{k}(t_0;t)$ is a finite union of smooth surfaces
(tori), and
\begin{equation}\label{lac-2-ax}
\int_{\widetilde S_k(t_0;t)}\nabla\Phi_k\cdot{\bf
n}\,dS=-\int_{\widetilde S_k(t_0;t)}|\nabla\Phi_k|\,dS<0,
\end{equation}
where $\n$ is the unit outward normal vector to
$\partial\widetilde W_k(t_0;t)$.

Note that $W_k(t_0;t)$ and
$S_k(t_0;t)$ are well defined for all $t\in(t_0,t_*)$ and $k\ge
 k_\circ= k_\circ(t_0)$. \  Now we are ready to prove the key estimate (which is analog of Lemma~3.8 from~\cite{kpr15annals}).

\bl
\label{ax-lkr11}{\sl
 Let $0<t_0<t_*$. Then there exists $k_*=k_*(t_0)$ such that
for every $k\ge k_*$ and for almost all $t\in(t_0,t_*)$
the inequality
\begin{equation}\label{mec}
\int\limits_{\widetilde S_k(t_0;t)}|\nabla\Phi_k|\,dS<\F t,
\end{equation}
holds with the constant $\F$ independent of  $t,t_0$,
and $k$. }
\el

\bpf Fix positive $t_0< t_*$, $\epsilon=\frac12 t_0$, and the corresponding parameters $\delta_\epsilon>0$, $k_\circ:=k_\epsilon\in\N$  such that
(\ref{fest1})--(\ref{fest2}) holds for all $k\ge k_\circ$ and  $h\in\HI\cap(0,\delta_\epsilon)$. Fix a number
${h_\circ}\in\HI\cap(0,\delta_\epsilon)$. Then for all $t\in(t_0,t_*)$ and all $k\ge k_\circ$ we can define the set~$S_k(t_0;t)$ as above. Moreover, for almost all $t\in(t_0,t_*)$ the set $S_k(t_0;t)$
 consists of finitely many
pairwise disjoint $C^1$-cycles (=$C^1$-smooth curves homeomorphic to the~circle) and $\Phi_k$ is differentiable (in classical sense)
at every point~$x\in S_k(t_0;t)$ with $\nabla\Phi_k(x)\ne0$.
The values $t\in(t_0,t_*)$ having the above property will be called
$k$-{\it regular}.

Respectively, for every regular value~$t\in(t_0,t_*)$ the set
$\widetilde S_{k}(t_0;t)$ is a finite union of smooth surfaces
(tori), and the inequality~(\ref{lac-2-ax}) holds.

The main idea of the proof of~(\ref{mec}) is quite simple: we will
integrate the equation
\begin{equation}\label{cle_lap**}\Delta\Phi_k=\omega_k^2+\frac1{\nu_k}\div(\Phi_k\ue_k)-
\frac{1}{\nu_k}\fe_k\cdot\ue_k\end{equation} over the suitable
domain $\Omega_{k}(t)$ with $\partial\Omega_{k}(t)\supset
\widetilde S_k(t_0;t)$ such that the corresponding boundary
integrals
\begin{equation}\label{r**1}
\biggr|\int_{\bigl(\partial\Omega_{k}(t)\bigr)\setminus \widetilde
S_k(t_0;t)}\nabla\Phi_k\cdot{\bf n}\,dS\biggr|
\end{equation}
\begin{equation}\label{r**2}
\frac1{\nu_k}\biggr|\int_{\bigl(\partial\Omega_{k}(t)\bigr)\setminus
\widetilde S_k(t_0;t)}\Phi_k\ue_k\cdot{\bf n}\,dS\biggr|
\end{equation}
are negligible. We split the construction of the domain
$\Omega_{k}(t)$ into two steps.

First of all, define the open set $$\DD_k(t):=W_k(t_0;t)\cup \overline{\DD_{h_\circ}}\setminus\Si.$$
Then by construction (see, e.g., (\ref{boundary00}), (\ref{boundary6})\,) we have
\begin{equation}\label{bound9}
\partial\DD_k(t)=\Si\cup S_k(t_0;t).
\end{equation}
Denote further
$$\Omega_k(t):=\wt\DD_k(t).$$
Then $\Omega_k(t)$ is the open set in the three dimensional space and
\begin{equation}\label{bound10}
\partial\Omega_k(t)=\Gamma\cup\wt{S_k}(t_0;t),
\end{equation}
where we denote $\Gamma:=\Gamma_0\cup\Gamma_1\cup\dots\cup\Gamma_M$.

\medskip
 By  direct calculations, (\ref{preeuler}) implies
$$\nabla\Phi_k=-\nu_k\curl\,\ov_k+\ue_k\times\ov_k+\fe_k=-\nu_k\curl\,\ov_k+\ue_k\times\ov_k+\lambda_k\nu_k^2\,\curl\,{\bf b}.$$

By the Stokes theorem, for any $C^1$-smooth closed surface
$S\subset\Omega$ and ${\bf g}\in W^{2,2}(\Omega)$ we have
$$\int_S\curl{\bf g}\cdot{\bf n}\,dS=0.$$  So, in particular,
\begin{equation}\label{uesk0}\int_S\nabla\Phi_k\cdot{\bf
n}\,dS=\int_S(\ue_k\times\ov_k)\cdot{\bf n}\,dS.
\end{equation}
 Recall that  by
the pressure normalization condition,
\begin{equation}\label{remo1}
\Phi|_{\Gamma}=0.\end{equation}

Our purpose on the next step is as follows: for arbitrary $\e>0$ and
for sufficiently large~$k$ to prove the estimates
\begin{equation}\label{remop1}
\biggr|\int_{\Gamma}\nabla\Phi_k\cdot{\bf
n}\,dS\biggr|=\biggl|\int_{\Gamma}(\ue_k\times\ov_k)\cdot{\bf
n}\,dS\biggr| <\e,
\end{equation}
\begin{equation}\label{remop2}
\frac1{\nu_k}\biggr|\int_{\Gamma}\Phi_k\ue_k\,dS\biggr|<\e.
\end{equation}

Recall that in our notation $\ue_k=\lambda_k\nu_k\h+\w_k$,
where $\w_k\in H(\Omega)$, $\|\w_k\|_{H(\Omega)}=1$, and $\h$ is a~solution to the Stokes problem with boundary value~$\mathbf a$ and forcing term~$\mathbf f$  (see~(\ref{f})\,). In
particular, we have
\begin{equation}\label{simp1}
\ue_k(x)\equiv\lambda_k\nu_k\h(x)\qquad\forall x\in\Gamma.
\end{equation}
To establish (\ref{remop2}),  we use the uniform boundedness
\begin{equation}\label{simp2}
\|\Phi_k\|_{L^3(\Omega)}+\|\nabla\Phi_k\|_{L^{3/2}(\Omega)}\le
C.
\end{equation}
From~(\ref{simp2}) and the weak
convergence $\Phi_k\rightharpoonup \Phi$ in
$W^{1,3/2}(\Omega)$ we easily have
\begin{equation}\label{remob0.4}
\Phi_k\to\Phi \qquad\mbox{ in \ }L^q(\Gamma)\ \ \ \forall
q\in[1,2).
\end{equation}
Thus by virtue of (\ref{remo1}),
\begin{equation}\label{simp3}
\int_{\Gamma}|\Phi_k|\,dS\to0\qquad\mbox{ as }k\to\infty.
\end{equation}Since \begin{equation}\label{est_i1} \|\h\|_{L^\infty(\Omega)}\le C
\|\h\|_{W^{2,2}(\Omega)}<\infty,
\end{equation}
by identity~(\ref{simp1}) we have
\begin{equation}\label{simp4}
\frac1{\nu_k}\biggr|\int_{\Gamma}\Phi_k\ue_k\,dS\biggr|=
\lambda_k\biggr|\int_{\Gamma}\Phi_k\h\,dS\biggr|\le
C\int_{\Gamma}|\Phi_k|\,dS \underset{k\to\infty}\to0,
\end{equation}
that implies the required estimate~(\ref{remop2}) for sufficiently
large~$k$.

To prove~(\ref{remop1}), we need also the uniform estimate
\begin{equation}\label{est_i2}
\|\nu_k\ue_k\|_{W^{2,3/2}(\Omega)}\le C,
\end{equation}
where $C$ is independent of~$k$ (this  inequality follows from the
construction, see~(\ref{stokes})\,).
Thus by Sobolev imbedding theorems
\begin{equation}\label{simp5}
\|\nu_k\nabla\ue_k\|_{L^1(\Gamma_{h})}\le C\qquad\forall
h\in[0,{h_\circ}]\ \ \forall k\in\N,
\end{equation}
where $C>0$ is independent of~$k,h$, here
$\Omega_h=\{x\in\Omega:\dist(x,\Gamma)\le h\}$,\ \
$\Gamma_h=\{x\in\Omega:\dist(x,\Gamma)=h\}$. Moreover, by
elementary calculations (\ref{est_i2}) implies the uniform
H\"{o}lder continuity of the function $[0,{h_\circ}]\ni h\mapsto
\|\nu_k\nabla\ue_k\|_{L^1(\Gamma_{h})}$, i.e., there exists a
constant $\sigma>0$ (independent of~$k$\,) such that
\begin{equation}\label{simp6}
\biggl|\int_{\Gamma_{h'}}\bigl|\nu_k\nabla\ue_k\bigr|\,dS-
\int_{\Gamma_{h''}}\bigl|\nu_k\nabla\ue_k\bigr|\,dS\biggr|\le
\sigma\bigl|h'-h''\bigr|^{\frac13}\qquad\forall
h',h''\in[0,{h_\circ}]\quad\forall k\in\N.
\end{equation}
From the last property and from the uniform boundedness of the
Dirichlet integral
\begin{equation}\label{simp7}
\|\nabla\ue_k\|_{L^{2}(\Omega)}\le 1+\lambda_k\nu_k\|\nabla\mathbf U\|_{L^{2}(\Omega)}\le 2
\end{equation}
(for sufficiently large~$k$\,) one can easily deduce that
\begin{equation}\label{simp8}
\sup\limits_{h\in[0,{h_\circ}]}\int_{\Gamma_{h}}\bigl|\nu_k\nabla\ue_k\bigr|\,dS\to0\qquad\mbox{
as }k\to\infty,
\end{equation}
in particular,
\begin{equation}\label{simp9}
\int_{\Gamma_{0}}\bigl|\nu_k\nabla\ue_k\bigr|\,dS\to0\qquad\mbox{
as }k\to\infty.
\end{equation}
 Then from the identity (\ref{simp1}) and the estimate (\ref{est_i1})  we have
\begin{equation}\label{simp10}
\biggl|\int_{\Gamma}(\ue_k\times\ov_k)\cdot{\bf n}\,dS\biggr|
=\lambda_k\nu_k\biggl|\int_{\Gamma}(\h\times\ov_k)\cdot{\bf
n}\,dS\biggr| \le
C\int_{\Gamma_{0}}\bigl|\nu_k\nabla\ue_k\bigr|\,dS\to0\qquad\mbox{
as }k\to\infty.
\end{equation}
Hence the required estimate~(\ref{remop1}) is proved.

\begin{center}
\includegraphics[height=6cm, width=8cm]{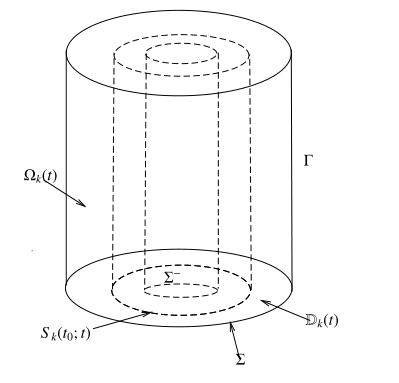}\\
{\small Figure 2: {\sl The domain $\Omega_k(t)$. }}
\end{center}

\medskip

Recall that by construction we have \ $\partial\Omega_k(t)=\Gamma\cup\wt{S_k}(t_0;t)$ \ (see~(\ref{bound10})\,). Integrating the
equation
\begin{equation}\label{cle_lap}\Delta\Phi_k=\omega_k^2+\frac1{\nu_k}\div(\Phi_k\ue_k)-\frac{1}{\nu_k}\fe_k\cdot\ue_k\end{equation}
over the domain $\Omega_{k}(t)$, we obtain
$$
\int_{ \widetilde S_k(t_0;t)}\nabla\Phi_k\cdot{\bf
n}\,dS+\int_{ \Gamma}\nabla\Phi_k\cdot{\bf
n}\,dS=\int_{\Omega_{k}(t)}\omega_k^2\,dx-\frac1{\nu_k}\int_{\Omega_{k}(t)}\fe_k\cdot\ue_k\,dx
$$
$$
+\frac{1}{\nu_k}\int_{ \widetilde
S_k(t_0;t)}\Phi_k\ue_k\cdot{\bf n}\,dS+\frac{1}{\nu_k}\int_{\Gamma}\Phi_k\ue_k\cdot{\bf n}\,dS
$$
\begin{equation}\label{cle_lap1}
=\int_{\Omega_{k}(t)}\omega_k^2\,dx-\frac1{\nu_k}\int_{\Omega_{k}(t)}\fe_k\cdot\ue_k\,dx+t{}\bar\F_k+\frac{1}{\nu_k}\int_{\Gamma}\Phi_k\ue_k\cdot{\bf n}\,dS,
\end{equation}
where $\bar\F_k=\lambda_k(\F_{0}+\dots+\F_M)$  and $\F_j=\int\limits_{\Gamma_j}\mathbf a\cdot\mathbf n\,dS$ is the corresponding boundary flux (here we use the
identity~$\Phi_k\equiv-t$ on $\widetilde S_k(t_0;t)$\,). In
view of (\ref{lac-2-ax}) and
(\ref{remop1})--(\ref{remop2})
 we can estimate
\begin{equation}\label{cle_lap2}
\int_{ \widetilde S_k(t_0;t)}|\nabla\Phi_k|\,dS \le
t\F_k+3\e+\frac1{\nu_k}\int_{\Omega_{k}(t)}\fe_k\cdot\ue_k\,dx-\int_{\Omega_{k}(t)}\omega_k^2\,dx
\end{equation} with $\F_k=|\bar\F_k|$. By definition, $\frac{1}{\nu_k}\|\fe_k\|_{L^{2}(\Omega)}=
\lambda_k{\nu_k}\|\fe\|_{L^{2}(\Omega)}\to 0$ as
$k\to\infty$. Therefore, using the uniform
estimate~$\|\ue_k\|_{L^6(\Omega)}\le\const$, we have
$$
\Big|\frac1{\nu_k}\int_{\Omega_{k}(t)}\fe_k\cdot\ue_k\,dx\Big|<\varepsilon
$$
for sufficiently large $k$. Then (\ref{cle_lap2}) yields
\begin{equation}\label{cle_lap3}
\int_{ \widetilde S_k(t_0;t)}|\nabla\Phi_k|\,dS <
t\F_k+4\e-\int_{\Omega_{k}(t)}\omega_k^2\,dx.
\end{equation}
Since $\F_k=\lambda_k|\F_{0}+\dots+\F_M|$ is uniformy bounded, and $\e$ could be taken arbitrary small, the last inequality implies the required estimate~(\ref{mec}) for sufficiently large~$k$.
\epf

\bl
\label{lem_Leray_fc} {\sl Assume that $\Om\subset\mathbb{R}^{2}\times\mathbb{T}_{\sigma}$ is a helical domain with $C^2$ boundary $\p\Om$, ${\bf f}={\rm curl}\,{\bf b}$, ${\bf b}\in W^{1,2}(\Om)$ and ${\bf a}\in W^{3/2,2}(\p\Om)$ are both helical functions, and $(\ref{ch})$ is valid. Then assumptions
${\bf(E\textendash NS\textendash H)}$ with \eqref{a} lead to a contradiction.}
\el

\bpf Fix positive $t_0< t_*$. Put $t_{i}=2^{-i}t_0$, \ $i=1,2,\dots$.
Take the  corresponding paprameters  $k_{*i}=k_{*}(t_i)$ such that
 the inequality~(\ref{mec}) holds for all $k\ge k_{*i}$ and for almost all $t\in(t_{i},t_*)$. In particular,
~(\ref{mec}) holds for all $k\ge k_{*i}$ and for almost all $t\in(t_{i},2t_i)$.

 For $k\ge k_{*i}$  put
$$
E_{ki}=\bigcup\limits_{t\in [t_{i},2t_i ]} \widetilde S_k(t_i;t).
$$
 By the Coarea formula (see, e.g, \cite{Maly}), for any integrable
 function $g:E_{ki}\to\R$ the equality
\begin{equation}\label{Coarea_Phi}\int\limits_{E_{ki}}g|\nabla\Phi_k|\,dx=
\int\limits_{t_i}^{2t_i}\int_{\widetilde
S_k(t_i;t)}g(x)\,d\Ha^2(x)\,dt
\end{equation}
holds. In particular, taking $g=|\nabla\Phi_k|$ and
using~(\ref{mec}), we obtain
\begin{equation}\label{Coarea_Phi2}
\int\limits_{E_{ki}}|\nabla\Phi_k|^2\,dx=
\int\limits_{t_i}^{2t_i}\int_{ \widetilde
S_{k}(t_i;t)}|\nabla\Phi_k|(x)\,d\Ha^2(x)\,dt
\le
\int\limits_{t_i}^{2t_i}\F t\,dt=
\frac\F2\bigl((2t_i)^2-(t_i)^2\bigr)\le 2\F t^2_i
\end{equation}
Now, taking  $g=1$ in (\ref{Coarea_Phi})  and using the H\"older
inequality we have
\begin{equation}\label{Coarea_Phi3}
\int\limits_{t_i}^{2t_i}\Ha^2\bigl( \widetilde
S_{k}(t_i;t)\bigr)\,dt= \int\limits_{E_{ki}}|\nabla\Phi_k|\,dx
\le
\biggl(\int\limits_{E_{ki}}|\nabla\Phi_k|^2\,dx\biggr)^{\frac12}
\bigl(\meas
(E_{ki})\bigr)^{\frac12}\le\sqrt{2\F t_i^2\meas
(E_{ki})}\le t_i\sqrt{2\F\meas
(E_{ki})}.
\end{equation}
 By construction (see the arguments after the formula~(\ref{boundary2})\,), for all $k$-regular values~$t$
the set $\widetilde
S_{k}(t_i;t)$ is a finite union of $C^1$-smooth  surfaces (tori) separating  $\Gamma=\Gamma_0\cup\dots\Gamma_M$
from
$\Gamma_{M+1}\cup\dots\cup\Gamma_N$. It implies, in particular,
that
$\Ha^2\bigl(\widetilde
S_{k}(t_i;t)\bigr)\ge C_*=C_*(\Omega)>0$. Therefore, by virtue of (\ref{Coarea_Phi3}) we have
\begin{equation}\label{Coarea_Phi3'}
C_*t_i\le t_i\sqrt{2\F\meas
(E_{ki})},
\end{equation}
in other words,
\begin{equation}\label{Coarea_Phi3''}
C_*\le \sqrt{2\F\meas
(E_{ki})}
\end{equation}
for all $i\in\N$ and for all $k\ge k_{*i}$.
By constructions,
 for all $k\ge k_{*i}$ the sets
$E_{ki}$, $E_{k(i-1)}$,  \dots, $E_{k2}$, $E_{k1}$ are pairwise disjoint.
Therefore, the measure of some $E_{ki}$ could be made arbitrary small (for sufficiently large~$i$ and $k$). This obviously contradicts the~estimate~(\ref{Coarea_Phi3''}). The Lemma is proved.
\epf

\subsubsection{If $(b)$ happens.}

Suppose now that the maximum of $\Phi$ is not attained on the boundary $\p\Om$\footnote{$\esssup\Phi=\oo$ is not excluded.}:
\be\label{bbb}
\max_{j=0,...,N}\check{p}_{j}<\esssup\limits_{x\in\Omega}\Phi(x).
\ee

Adding  a constant to the pressure, we can assume without loss of generality that
\be\label{bbb1}
\max_{j=0,...,N}\check{p}_{j}<p_*<0<\esssup\limits_{x\in\Omega}\Phi(x).
\ee

The proof for this case can be carried out with the same arguments as in the previous subsection with obvious simplifications. Let us describe some details. We start from  the following simple fact.

\bl
\label{lem_PR} {\sl Under assumptions~$(\ref{bbb1})$ there exists a~straight segment~$F\subset \DD$ such that
$F\cap A_{\we}=\emptyset$ \,(i.e., $F$ consists of the regular points of~$\Phi$ and  $\we$, see~Lemma~\ref{sobolev}\,), and
$$0<\inf\limits_{x\in F}\Phi(x),$$
moreover, the uniform convergence
\begin{equation}\label{ucc3}\Phi_{k}|_F\rightrightarrows \Phi|_F
\end{equation}
holds.}
\el

This fact follows easily from the definition of Sobolev spaces and from the weak convergence~$\Phi_k\rightharpoonup\Phi$ in $W^{1,3/2}(\DD)$ (see, e.g., the proof of Theorem~3.2 in \cite{aiumj84} for details), so we omit its proof here.

Fix the segment   $F$ from Lemma~\ref{lem_PR}. Denote $t_*=-p_*$ and fix a~positive $t_0<t_*$.
Using the arguments from the~previous subsection, we can find a~small parameter $h_\circ>0$ and a~number $k_\circ\in\N$ such that
\begin{equation}\label{boundary-bb2}
\forall t\in(t_0,t_*)\ \ \forall k\ge k_\circ\quad
\biggl(\Phi_k|_F>0,\quad \Phi_k|_{\Si^-_{h_\circ}}< -t\biggr),
\end{equation}
where now by definition $\Si^-_{h}=\{x\in\DD:\dist(x,\partial\DD)=h\}$ (see Fig.~3).

\begin{center}
\includegraphics[height=6cm, width=7cm]{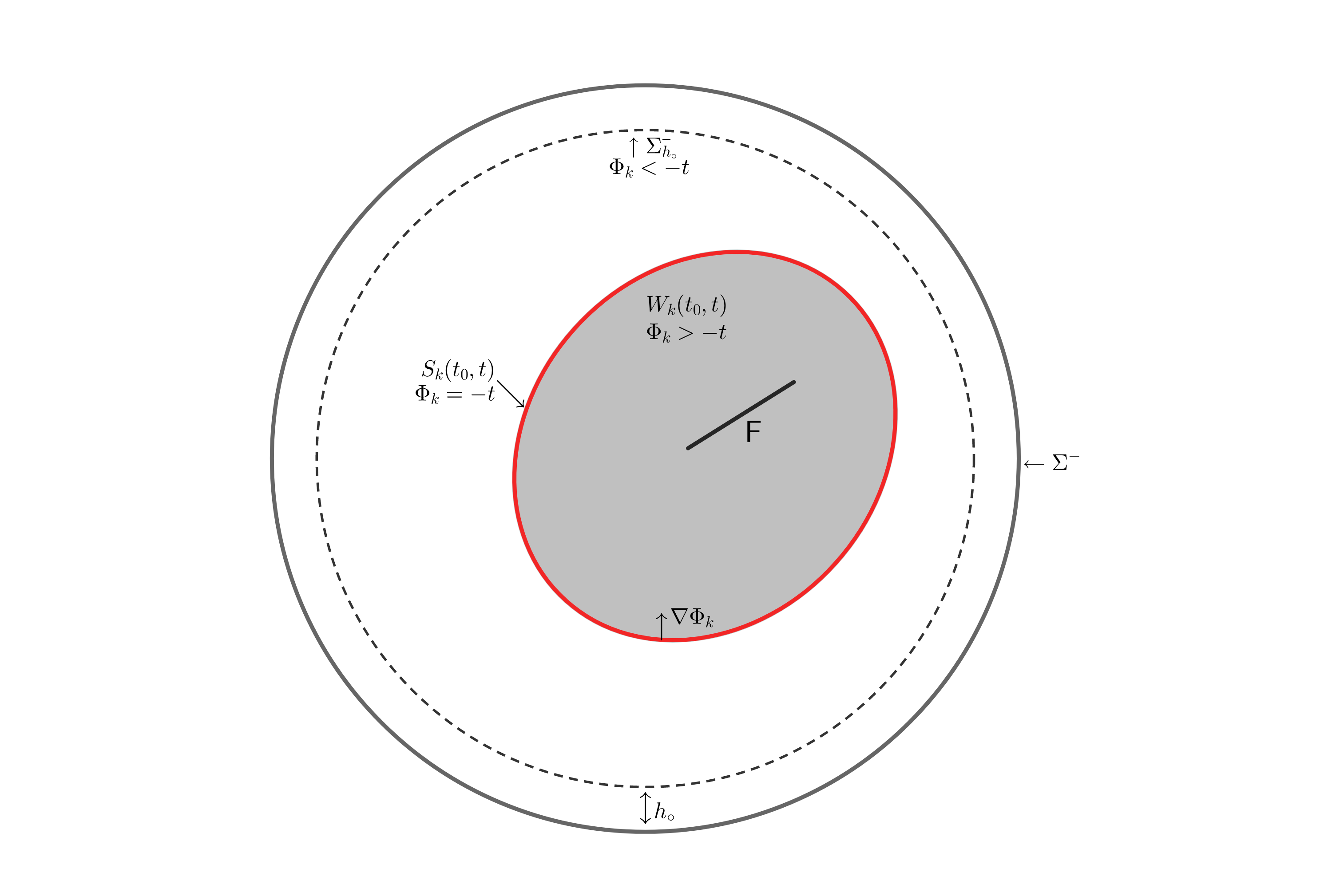}\\
{\small Figure 3. }
\end{center}

Further, using the same arguments, for almost all $t\in(t_0,t_*)$ we can find a set~$S_k(t_0;t)$
consisting of finite disjoint family of $C^1$ smooth closed curves (cycles) such that $\Phi_k\equiv -t$ on
$S_k(t_0;t)$. Moreover, there is an open set $W_k(t_0;t)\subset\DD$  satisfying the relations
\begin{equation}\label{boundary-bbb1}W_k(t_0;t)\supset  F, \end{equation}
\begin{equation}\label{boundary-bbb11}\partial W_k(t_0;t)=S_k(t_0;t)\end{equation}
(cf. with~(\ref{boundary6})\,), and
\begin{equation}\label{bbb-lac-2-ax}
\int_{\widetilde S_k(t_0;t)}\nabla\Phi_k\cdot{\bf
n}\,dS=-\int_{\widetilde S_k(t_0;t)}|\nabla\Phi_k|\,dS<0,
\end{equation}
where $\n$ is the unit outward normal vector to
$\partial\widetilde W_k(t_0;t)$. Here, by construction, $\widetilde S_{k}(t_0;t)$,  is a finite union of smooth disjoint  surfaces (tori).

Further we can prove the estimate~(\ref{mec}) for our case integrating the same identity~(\ref{cle_lap**}) over the domain~$\Omega_k(t)=\widetilde W_k(t_0;t)$ with $\partial\Omega_k(t)=\widetilde S_k(t_0;t)$. Note that now the proof is even much simpler since we have no boundary integrals over subsets
of~$\partial\Omega$~--- in other words, now we do not need to prove estimates of type~(\ref{r**1})--(\ref{r**2}) and (\ref{remop1})--(\ref{remop2})

After the estimate~(\ref{mec}) is obtained, we can derive a~contradiction  in exactly the same way as in Lemma~\ref{lem_Leray_fc} of the~previous section.

\bigskip Now we can summarize the results of the last two subsections in the following statement.

\bl
\label{lem_Leray_fc_3-bbb} {\sl Assume that $\Om\subset\mathbb{R}^{2}\times\mathbb{T}_{\sigma}$ is a helical domain with $C^2$ boundary $\p\Om$, ${\bf f}={\rm curl}\,{\bf b}$, ${\bf b}\in W^{1,2}(\Om)$ and ${\bf a}\in W^{3/2,2}(\p\Om)$ are both helical functions, and $(\ref{ch})$ is valid.
 Let ${\bf(E\textendash NS\textendash H)}$ be fulfilled. Then each assumptions
$(\ref{a})$ or $(\ref{b})$ lead to a contradiction.
}
\el

\bigskip {\bf Proof of Theorem \ref{main}.} Let the hypotheses
of Theorem~\ref{main} be satisfied. Suppose that its assertion
fails. Then, by Lemma~\ref{NSkk2}, there exist $ \ve, p$
and a sequence $(\ue_k,p_k)$ satisfying~{\bf(E\textendash NS\textendash H)}, and by Lemma~\ref{lem_Leray_fc_3-bbb}
these
assumptions lead to a contradiction.    \qed


\end{document}